\input amstex 
\documentstyle{amsppt}
\input bull-ppt
\keyedby{bull495/car}

\topmatter
\cvol{31}
\cvolyear{1994}
\cmonth{July}
\cyear{1994}
\cvolno{1}
\cpgs{44-49}
\ratitle
\title Lyapunov Theorems for Banach Spaces\endtitle
\author Y. Latushkin and S. Montgomery-Smith\endauthor
\address Department of Mathematics, University of 
Missouri, Columbia, Missouri
65211\endaddress
\ml\nofrills {\it E-mail address}{\rm, {\smc Y. 
Latushkin}:}\enspace
\tt mathyl\@mizzou1.missouri.edu
\mlpar
{\it E-mail address}{\rm, {\smc S. 
Montgomery-Smith}:}\enspace
stephen\@mont.cs.missouri.edu\endml
\date February 19, 1993\enddate
\subjclass Primary 4706, 47B38; Secondary 34D20, 
34G10\endsubjclass
\keywords Hyperbolicity, evolution family, exponential 
dichotomy, weighted
composition operators, spectral mapping theorem\endkeywords
\abstract We present a spectral mapping theorem for 
semigroups on any Banach
space $E$. From this, we obtain a characterization of 
exponential dichotomy for
nonautonomous differential equations for $E$-valued 
functions. This 
characterization is given in terms of the spectrum of the 
generator of the
semigroup of evolutionary operators.\endabstract
\endtopmatter

\document

\heading 1. Introduction\endheading

Let us consider an autonomous differential equation 
$y'=Ay$ in a Banach space
$E$, where $A$ is a generator of continuous semigroup 
$\{e^{tA}\}_{t\geq0}$;
that is, the solution of the differential equation 
satisfies $y(t)=e^{tA}y(0)$,
$t\geq0$. A classical result of A. M. Lyapunov (see, e.g., 
\cite{DK}) shows
that for {\it bounded\/} $A$, the spectrum $\sigma(A)$ of 
$A$ is responsible
for the asymptotic behavior of $y(t)$. For example, if 
$\operatorname{ Re}
\sigma(A)<0$, then the trivial solution is uniformly 
asymptotically stable,
that is, $\|e^{tA}\|\rightarrow 0$ as $t\rightarrow \infty 
$. This fact follows
from the spectral mapping theorem (see, e.g., \cite{N, p.\ 
82}):

$$\sigma(e^{tA})\backslash \{0\}=\exp(t\sigma(A)),\quad 
\quad t\neq0,\tag1$$
which always holds for bounded $A$.

For unbounded $A$, equation (1) is not always true. 
Moreover, there are
examples of generators $A$ (see \cite{N, p.\ 61}) such 
that even $
\operatorname{ Re}(\sigma(A))\leq s_0<0$ does not 
guarantee $\sigma(e^A)
\subset  \Bbb D=\{|z|<1\}$ or $\|e^{tA}\|\rightarrow 0$ as 
$t\rightarrow
\infty $. Since $\sigma(A)$ does not characterize the 
asymptotic behavior of
the solutions $y(t)$, we would like to find some other 
characterization
that still does not involve solving the differential 
equation (i.e., finding $
\sigma(e^{tA}))$.

In this article we solve precisely this problem in the 
following manner.
Consider the space $L_p(\Bbb R;E)$ of $E$-valued functions 
for $1\leq p<\infty
$ and the semigroup $\{e^{tB}\}_{t\geq0}$ of evolutionary 
operations (also
called weighted translation operators)
$$(e^{tB}f\,)(x)=e^{tA}f(x-t),\quad \quad t\geq0,\tag2$$
generated by the operator $B=-\frac d{dx}+A$, $x\in \Bbb 
R$. It turns out that
it is $\sigma(B)$ in $L_p(\Bbb R;E)$ that is responsible 
for the asymptotic
behavior of $y(t)$ in $E$. For example, $\operatorname{ 
Re}(\sigma(B))<0$ on
$L_p(\Bbb R;E)$ implies that $\|e^{tA}\|\rightarrow 0$ as 
$t\rightarrow \infty
$ on $E$.

We will also consider the well-posed nonautonomous 
equation $y'=A(t)y(t)$.
Instead of the semigroup given by (2), we consider in 
$L_p(\Bbb R;E)$ the
semigroup
$$(e^{tD}f\,)(x)=U(x,x-t)f(x-t),\qquad x\in \Bbb R,\ 
t\geq0\.\tag3$$
Here $U(t,s)$, $t\geq s$, is the evolutionary family 
(propagator) for the
nonautonomous equation. We will show that $\sigma(D)$ 
characterizes the 
asymptotic behavior of $y(t)$.

The  order of proofs will be as follows. We will first 
show the spectral
mapping theorem for the autonomous case (2). We will also 
consider a similar 
theorem for the evolutionary semigroup in the space of the 
periodic functions.
This theorem will give us a variant of Greiner's spectral 
mapping theorem (see
\cite{N, p.\ 94}) for any $C_0$-semigroup $\{e^{tA}\}$ in 
a Banach space. This
variant also is a direct generalization of Gerhard's 
spectral mapping theorem
in Hilbert space for generators with resolvent bounded 
along $i\Bbb R$ (see 
\cite{N, p.\ 95}). Then we will obtain the spectral 
mapping theorem for the
nonautonomous case (3) using a simple change of variables 
argument to reduce it
to the autonomous case (2).

We will be considering not only stability but also the 
exponential dichotomy
(hyperbolicity) for the solutions of the equation 
$y'=A(t)y(t)$ in $E$. In the 
theory of differential equations with bounded 
coefficients, exponential
dichotomy is an important tool used, for example, in 
proving instability 
theorems for nonlinear equations and determining existence 
and uniqueness of
bounded solutions and Green's functions (see, e.g., 
\cite{DK}). The spectral
mapping theorem given here for the semigroup (3) allows 
one to extend these
ideas to the case of unbounded coefficients.

In turns out that the condition $0\notin \sigma(D)$---or 
equivalently
$\sigma(e^{tD})\cap \Bbb T=\varnothing$, $t>0$, $\Bbb 
T=\{|z|=1\}$ on
$L_p(\Bbb R;E)$---is equivalent to the hyperbolic behavior 
of a special kind
for the solutions. We will call this {\it spectral 
hyperbolicity\/}. Note that
if the $U(t,s)$ are invertible for $(t,s)\in \Bbb R^2$, 
that is, (3) is
extendible to a group, then the spectral hyperbolicity is 
the same as the
exponential dichotomy (or hyperbolicity) in the usual 
(see, e.g., \cite{DK})
sense. Therefore, the spectrum $\sigma(e^{tD})$ for 
nonperiodic $A(\boldcdot)$
plays the same role in the description of exponential 
dichotomy as the spectrum
of the monodromy operator does in the usual Floquet theory 
for the periodic
case. However, ordinary hyperbolicity is not equivalent 
(see \cite R) to spectral
hyperbolicity in the {\it semi\/}group case and thus 
cannot be characterized in
terms of $\sigma(e^{tD})$ only. For a Hilbert space, we 
were able to
characterize hyperbolicity in terms of other spectral 
properties of $e^{tD}$.

Finally, the results of this article can be generalized to 
the case of the
variational equation $y'(t)=A(\varphi^tx)y(t)$ for a flow 
$\{\varphi^t\}$ on
a compact metric space $X$ or for a linear skew-product flow
$\hat{\varphi}^t\:X\times E\rightarrow X\times
E\:(x,y)\mapsto(\varphi^tx,\Phi(x,t)y)$, $t\geq0$ (see 
\cite{CS, H, LS, SS} and
references contained therein). Here $\Phi\:X\times \Bbb R_+
\rightarrow L(E)$ is
a cocycle over $\varphi^t$, that is, $\Phi(x,t+s)=
\Phi(\varphi^tx,s)\Phi(x,t)$,
$\Phi(x,0)=I$. Let us recall (see \cite{SS}) that one of 
the purposes of the 
theory of linear skew-product flows was to aid in studying 
the equation
$y'=A(t)y$ for the case of almost periodic $A(\boldcdot)$. 
To answer the
question when $\hat \varphi^t$ is hyperbolic (or Anosov), 
instead of (3) one
considers the semigroup of so-called weighted composition 
operators (see
\cite{CS, J, LS}) on $L_p(X;\mu;E)$:
$$(T^tf\,)(x)=\left(\frac{d\mu\circ 
\varphi^{-t}}{d\mu}\right)^{1/p} \Phi(
\varphi^{-t}x,t)f(\varphi^{-t}x),\quad \quad x\in X,\ 
t\geq0\.\tag4$$
Here $\mu$ is a $\varphi^t$-quasi-invariant Borel measure 
on $X$. As above, the 
condition $\sigma(T^t)\cap\Bbb T=\varnothing$ is 
equivalent to the spectral
hyperbolicity of the linear skew-product flow $\hat 
\varphi^t$. The spectral
hyperbolicity coincides with the usual hyperbolicity if 
$\Phi(x,t)$, $x\in X$,
$t\geq0$, are invertible or compact operators. Unlike the 
finite-dimensional case
(see \cite M), the hyperbolicity of $\hat \varphi^t$ does 
not generally imply
the condition $\sigma(T^t)\cap \Bbb T=\varnothing$. For 
Hilbert space the 
condition of hyperbolicity of $\hat \varphi^t$ can also be 
described in terms
of other spectral properties of $T^t$. We will not include 
these
generalizations in this paper.

We point out that the investigation of evolutionary 
operators (2) and (3) has a
long history \cite{Ho}. 
More recently, significant progress has been made in 
\cite{BG, N,
P, R} (see \cite R for detailed bibliography), and this 
list does not pretend
to be complete. A detailed investigation of weighted 
composition operators (4)
for Hilbert spaces $E$ and connections with the spectral 
theory of linear
skew-product flows \cite{SS} and other questions of 
dynamical systems theory
and a bibliography may be found in \cite{LS}.

\rem{Notation} $L(E)$ (correspondingly $L_s(E))$ denotes 
the set of bounded
operators on $E$ with the uniform (correspondingly strong) 
topology; $\rho(A)$
denotes the resolvent set of the operator $A;\ |$ denotes 
the restriction of an 
operator; $C_b(\Bbb R;E)$ denotes the space of continuous 
bounded $E$-valued 
functions on $\Bbb R$ with the supremum norm, and 
$C_b^0(\Bbb R;E)$ denotes the
subspace of functions vanishing at infinity.\endrem

\heading 2. Autonomous case\endheading

Let $A$ be a generator of any $C_0$-semigroup in a Banach 
space $E$. Consider
the $C_0$-semigroup (2) on the space $L_p([0,2\pi);E)$, 
$1\leq p<\infty $,
that is, $(e^{tB}f\,)(x)=e^{tA}f((x-t)(\operatorname{ 
mod}2\pi))$.
\proclaim{Theorem 1} The following are equivalent\RM:

\roster
\item \<$1\in \rho(e^{2\pi A})$ on $E$\RM;
\item \<$0\in \rho(B)$ on $L_p([0,2\pi);E)$\RM;
\item \<$1\in \rho(e^{2\pi B})$ on $L_p([0,2\pi);E)$.
\endroster
\endproclaim

In the main part of the proof (2)$\Rightarrow$(1), we 
modify the idea of
\cite{CS}. Let us assume that (2) is fulfilled, but for 
each $0<\varepsilon <
\frac12$, there is a $y\in E$ such that $\|e^{2\pi 
A}y-y\|<\varepsilon $ and
$\|y\|=1$ (and hence $\|e^{2\pi A}y\|\geq \frac12)$. Let 
$\rho(x)=\frac
3{2\pi}x-1$ for $x\in [
\frac{2\pi}3,\frac{4\pi}3)$, $\rho(x)=0$ for 
$x\in[0,\frac{2\pi}3)$, and
$\rho(x)=1$ for $x\in [\frac{4\pi}3,2\pi)$. Define the 
function $f\in L_p([0,2
\pi);E)$ by $f(x)=(1-\rho(x))e^{(2\pi+x)A}y+
\rho(x)e^{xA}y$, $x\in[0,2\pi)$.
Then for $c=\max\{\|e^{xA}\|\:x\in [0,2\pi)\}$, one has 
$\|e^{2\pi
A}y\|=\|e^{(2\pi-x)A}e^{xA}y\|\leq c\|e^{xA}y\|$. But, in 
contradiction with
(2), 
$$\|f\,\| ^p_p\geq 
\int^{2\pi}_{4\pi/3} \|e^{xA}y\|^p\,dx\geq \frac{2\pi 
}3c^{-p}\|e^{2\pi
A}y\|^p\geq \frac{2\pi }3c^{-p}2^{-p}$$
and
$$\|Bf\,\|_p\leq \frac{2\pi }3c\varepsilon \.$$

Theorem 1 implies the following variant of Greiner's 
spectral mapping theorem
(see \cite{N, p.\ 94}) for a $C_0$-semigroup $\{e^{tA}\}$ 
in Banach space $E$.

\proclaim{Theorem 2} $1\in \rho(e^{2\pi A})$ if and only 
if \RM{(a)} $i\Bbb Z
\subset  \rho(A)$ and \RM{(b)} there is a constant $C$ 
such that for any finite 
sequence $\{y_k\}\subset  E$
$$\left\|\sum_{ 
k}(A-ik)^{-1}y_ke^{-ikx}\right\|_{L_p([0,2\pi);E)}\leq
C\left\|\sum_{ k}y_ke^{-ikx}
\right\|_{L_p([0,2\pi);E)}\.$$
\endproclaim

Obviously, if $E$ is a Hilbert space and $p=2$, then 
Parseval's identity allows
one to replace (b) by the condition $\sup\{\|(A-ik)^{-1}\|:k
\in \Bbb Z\}<\infty $. This gives the famous spectral 
mapping theorem of
Gerhard \cite{N, p.\ 95}.

Let us consider (2) on $L_p(\Bbb R;E)$, $1\leq p<\infty $. 
The spectrum
$\sigma(B)$ now is invariant under the translations along 
the imaginary axis.
Moreover, we have the following result.

\proclaim{Theorem 3} For each $t>0$ the following are 
equivalent\RM:
\roster
\item \<$\sigma(e^{tA})\cap \Bbb T=\varnothing$ on $E$\RM;
\item \<$0\in \rho(B)$ on $L_p(\Bbb R;E)$\RM;
\item \<$\sigma(e^{tB})\cap \Bbb T=\varnothing$ on 
$L_p(\Bbb R;E)$.
\endroster\endproclaim

Thus, the spectral mapping theorem is valid for (2). If 
$\{e^{tA}\}_{t\in 
\Bbb R}$ is a group, then (1) is the same as the 
exponential dichotomy of the
autonomous equation $y'=Ay$ on $\Bbb R$.

\heading 3. Nonautonomous case\endheading

Consider the well-posed nonautonomous equation 
$y'(t)=A(t)y(t)$. By
``well-posed'' we mean that we assume the existence of a 
jointly strongly
continuous evolutionary family $U(t,s)\in L_s(E)$, $t\geq 
s$, with the 
properties $U(t,t)=I$, $U(t,r)=U(t,s)U(s,r)$, and 
$\|U(r,s)\|\leq
Ce^{\beta(t-s)}$, $t\geq r\geq s$. In fact, $U$ is a 
propagator for the 
equation $y'(t)=A(t)y(t)$, that is, $y(t)=U(t,s)y(s)$. The 
spectral mapping
theorem is valid for (3).

\proclaim{Theorem 4} 
Let \RM{(3)} be a $C_0$-semigroup on $L_p(\Bbb R;E)$, $1\leq
p<\infty $. Then $\sigma(D)$ is invariant under 
translations along the imaginary
axis and the following are equivalent\RM:
\roster
\item \<$0\in \rho (D)$ on $L_p(\Bbb R;E)$\RM;
\item \<$\sigma(e^{tD})\cap \Bbb T=\varnothing$ on 
$L_p(\Bbb R;E)$, $t>0$.
\endroster
\endproclaim

To outline the proof of (1)$\Rightarrow$(2), let us 
consider the semigroup
$(e^{tB}h)(s,x)=U(x,x-t)h(s-t,x-t)$, $(s,x)\in \Bbb R^2$, 
$t>0$, on the space
$L_p(\Bbb R\times \Bbb R;E)=L_p(\Bbb R;L_p(\Bbb R;E))$ and 
perform a change of
variables $u=s+x$, $v=x$. More precisely, consider the 
isometry $J$ on $L_p(
\Bbb R \times \Bbb R;E)$ defined by $(Jh)(s,x)=h(s+x,x)$. 
Then one has that
$Je^{tB}=(I\otimes e^{tD})J$  and $JB=(I\otimes D)J$, 
where $A_1\otimes A_2$
means that $A_1$ acts on $h(\boldcdot, x)$ and $A_2$ acts 
on $h(s,\boldcdot)$.
Since $e^{tB}$ can be written as 
$(e^{tB}f\,)(s)=e^{tD}f(s-t)$ for $f\:\Bbb R 
\rightarrow L_p(\Bbb R;E)\:s\mapsto h(s,\boldcdot)$, one 
can apply  to $e^{{tB}}$
the part  (2)$\Rightarrow$(3) of Theorem 3.

\dfn{Definition} The evolutionary family 
$\{U(x,s)\}_{x\geq s}$ is called {\it
hyperbolic\/} if there exists a projection-valued function 
$P\:\Bbb R 
\rightarrow L(E)$, $P\in C_b(\Bbb R;L_s(E))$, and $M$, 
$\lambda>0$ such that for all
$x\geq s$:
\roster 
\item \<$P(x)U(x,s)=U(x,s)P(s)$; and
\item \<$\|U(x,s)y\|\leq Me^{-\lambda(x-s)}\|y\|$ if $y\in 
\operatorname{
Im}P(s)$;
\item"\phantom{(2)}" \<$\|U(x,s)y\|\geq 
M^{-1}e^{\lambda(x-s)}\|y\|$ if
$y\in\Ker P(s)$.
\endroster
The evolutionary family $\{U(x,s)\}_{x\geq s}$ is called 
{\it spectrally
hyperbolic\/} if, in addition:
\roster
\item"(3)" \<$\operatorname{ Im}U(x,s)|\Ker P(s)$ is dense 
in $\Ker P(x)$.
\endroster\enddfn

Note that the second inequality in (2) implies only left 
invertibility of the
restriction $U(x,s)|\Ker P(s)$, while (3) guarantees its 
invertibility. If the
evolutionary family $\{U(x,s)\}_{(x,s)\in \Bbb R^2}$ 
consists of invertible
operators, then spectral hyperbolicity is equivalent to 
the hyperbolicity and is
the same as exponential dichotomy \cite{DK} of the 
equation $y'=A(t)y$ on $
\Bbb R$. If $\dim \Ker P(s)<d<\infty $, then obviously (2) 
always implies (3).
This also happens if the $U(x,s)$ are compact operators in 
$E$ \cite{R1}.
\proclaim{Theorem 5} The evolutionary family $\{U(x,s)\}$ 
on the separable 
Banach space $E$ is spectrally hyperbolic if and only if 
$0\in \rho(D)$ on
$L_p(\Bbb R;E)$.\endproclaim

\rem{Remark} The space $L_p$ may be replaced by the space 
$C([0,2\pi);E)$ in 
Theorem 1 and by $C^0_b(\Bbb R;E)$ in Theorems 2 to 5. The 
separability 
assumption in Theorem 5 was recently removed \cite{LR}.
\endrem

A remarkable observation in \cite R shows that the 
hyperbolicity of
$\{U(x,s)\}_{x\geq s}$ (unlike the spectral hyperbolicity) 
does not generally
imply (2) in Theorem 4 for infinite-dimensional $E$. 
However, we are able to
give the following characterization of hyperbolicity under 
the assumptions that
$E$ is a Hilbert space and that for some $r>0$ the 
function $x\mapsto U(x+r,x)$
is a continuous function from $\Bbb R$ to $L(E)$.

\proclaim{Theorem 6} The evolutionary family 
$\{U(x,s)\}_{x\geq s}$ is
hyperbolic on a separable Hilbert space $E$ if and only if 
there exists a
projection $\scr P$ on $L_2(\Bbb R;E)$ such that for some 
$t>0$\RM:
\roster
\item \<$e^{tD}\scr P=\scr P e^{tD}$\RM;
\item\<$\sigma(e^{tD}|\operatorname{ Im}\scr P)\subset  
\Bbb D$\RM;
\item \<$e^{tD}|\Ker \scr P$ is left invertible and 
$\sigma((e^{tD}|\Ker \scr
P)^\dagger)\subset  \Bbb D$\RM;
\item \<$\Ker \scr P\ominus \bigcap_{n\geq 
0}\operatorname{ Im}(e^{ntD}|\Ker 
\scr P)$ is invariant with respect to multiplications by 
the functions from
$C_b(\Bbb R;\Bbb R)$.
\endroster
Each projection $\scr P$ with these properties has a form 
$(\scr
Pf)(x)=P(x)f(x)$ for a projection-valued function $P\in 
C_b(\Bbb R;L(E))$.
\endproclaim

For the left-invertible operator $T$ the notation 
$T^\dagger$ stands for its
Moore-Penrose left inverse: $T^\dagger u=v$ if $u=Tv$, and 
$T^\dagger u=0$ for
$u\bot \operatorname{ Im}T$. Note that (1), (2), (3) imply 
the left
invertibility of $zI-e^{tD}$ for all $z\in \Bbb T$ and the 
formula $\scr P=
\frac 1{(2\pi i)}\int_{\Bbb T}(zI-e^{tD})^\dagger
\,dz$ which  gives the Riesz projection on 
$\sigma(e^{tD})\cap \Bbb D$ if $
\sigma(e^{tD})\cap \Bbb T=\varnothing$.

\enddocument